\documentclass[12pt]{article}
\usepackage{amssymb} 
\usepackage{amsmath}
\newtheorem{teo}{Theorem}

\newtheorem{rem}{Remark}
\newtheorem{lem}{Lemma}

\title{ Wavelet frames with matched masks
}
\author{ Elena A. Lebedeva
\footnote{Mathematics and Mechanics Faculty, Saint Petersburg State University,
Universitetsky prospekt, 28, Peterhof,  Saint Petersburg, 198504, Russia
 }
}
\date{
 ealebedeva2004@gmail.com
}

\begin{document}
\maketitle

\newcommand{\nul}{{\bf0}}
\newcommand{\rd}{{\mathbb R}^d}
\newcommand{\zd}{{\mathbb Z}^{d}}
\renewcommand{\r}{{\mathbb R}}
\newcommand{\z} {{\mathbb Z}}
\newcommand{\cn} {{\mathbb C}}
\newcommand{\n} {{\mathbb N}}

\begin{abstract} In the paper we design a Parseval wavelet frame  with a compact support and many vanishing moments.  The corresponding refinement mask approximates an arbitrary continuous periodic function $f$, $f(0)=1$. The refinable function has stable integer shifts. 
              
\end{abstract}

\textbf{Keywords} refinement mask, unitary extension principle, Parseval wavelet frame. 

\textbf{AMS Subject Classification}:  42C40

Wavelet frame as an overcomplete system is a very popular tool for a representation of a function. Especially, it is convenient to deal with a Parseval frame (a tight frame with a frame bound $A=1$), because the canonical dual frame coincides with the frame itself. The goal of our paper is to design a Parseval wavelet frame  with the following  specific approximation property: the refinement mask approximates an arbitrary continuous periodic function in some sense. For brevity we refer the mask with  this property as ``matched mask'' at the title of the article. In other words, we study density of the set of all refinement masks of framelets. Moreover, our wavelet system has other popular, important approximation properties. It has a compact support, many vanishing moments and the corresponding refinable function has stable integer shifts. 
To state the problem we are motivated by works due to Larson \cite{Larson}, Bownik \cite{Bownik}, Cabrelly, Molter \cite{CM}. In  \cite{Larson} the following question is stated: is the set of all wavelet generators dense in $L_2(\mathbb{R})$? In \cite{Bownik} and \cite{CM} a affirmative answer is given for framelets. 
In the paper we ask the same question for refinement masks.  

We recall briefly necessary notions and statements and stress on the corresponding properties of refinement masks. 
  To design frames we apply the unitary extension principle by Ron and Shen \cite{RS}. We recall it here. In the sequel, we denote by $C$ and $L_{p}$ the spaces of $2\pi$-periodic continuous functions and $2\pi$-periodic  functions $f$ such that $|f|^p$ is integrable respectively.  $\mathbb{T} = [-\pi,\,\pi)$   denotes the one-dimensional torus.

\textbf{The unitary  extension principle   }(UEP): 
{\it 
	Let $\varphi \in L_2(\mathbb{R})$ (the refinable function) and assume that there exists a function $m_0 \in L_{2}$  (the refinement mask) such that 
	\begin{equation}
	\label{refeq}
	\widehat{\varphi}(2 \xi) = m_0(\xi) \widehat{\varphi}(\xi), 
	\end{equation}
	and $\lim_{\xi \to 0}\widehat{\varphi}(\xi) =1.$ 
	Let $m_1,\dots, m_q \in L_{2}$ (the wavelet masks) and define $\psi_1,\dots,\psi_q\in L_2(\mathbb{R})$ (the wavelet generators) by 
	\begin{equation}
		\label{psi}
	\widehat{\psi_r}(2 \xi) = m_r(\xi) \widehat{\varphi}(\xi), \ \ r=1,\dots,q.
\end{equation}
	If for a.e. $\xi \in \mathbb{T}$
	\begin{equation}
		\label{matr}
		\left\{
	\begin{array}{l}
		\sum\limits_{r=0}^{q} |m_r(\xi)|^2 = 1, \\
		\sum\limits_{r=0}^{q} m_r(\xi)\overline{m_r(\xi+\pi)} = 0,
	\end{array}
	\right.
\end{equation}
	then the functions $\displaystyle \left\{\psi_{r,j,k}\right\}_{j,k\in \mathbb{Z},r=1,\dots,q}$ form a Parseval frame for $L_2(\mathbb{R}).$ Here we use the traditional notation $\psi_{r,j,k}(x):=2^{j/2}\psi_r(2^j x+k).$ 
	}

We intend to provide approximation of $L_2(\mathbb{R})$ by the multiresolution analysis (MRA) $\{V_j\}_{j\in\mathbb{Z}}$ and by  the wavelet frame system as well. 
We recall \cite{DHRS} that the MRA  provides approximation
order $N_{\varphi}$, if for every $f\in W^{N_{\varphi}}_2(\mathbb{R})$ 
$$
{\rm dist}(f,\,V_j) = O(2^{-j N_{\varphi}}).
$$
Let $\displaystyle Q_j(f) = \sum_{i< j}\sum_r\sum_k  \left\langle f,\,\psi_{r,i,k}\right\rangle \psi_{r,i,k}$
The wavelet system $\psi_{r,j,k}$ provides approximation order $N_{\psi}$, if
$$
\|f-Q_jf\|=O(2^{-j N_{\psi}}).
$$
It is clear that $N_{\varphi}\ge N_{\psi}.$ 
The approximation orders $N_{\varphi}$ and $N_{\psi}$ are controlled by the order of vanishing moments. To be more  precise, we recall that the wavelet system $\{\psi_{r,j,k}\}$ has vanishing moments of order $N_0$, if $\widehat{\psi_r}^{(k)}(0)=0$ for $k=0,\dots,N_0,$ $r=1,\dots,q.$  In \cite{DHRS}, it is proved that 
$$
N_{\varphi}\ge N_0 \ \  \mbox{ and } \ \ N_{\psi} = \min\{N_{\varphi},\,N\},
$$
where 
$N$ is the order of the zero of $1-|m_0|^2$   at the origin. The first equation in (\ref{matr}) and (\ref{psi})  
implies that $N=2N_0.$ Therefore,
$
N_0\leq N_{\psi}\leq 2 N_0.
$
  Therefore,  the order of vanishing moments is completely defined by the behavior of the refinement mask $m_0$ at the origin.  
  
	A refinable function $\varphi$ might not generate an MRA. That is why we are interested in the case when the refinable function $\varphi$ has stable integer shifts (see \cite[Definition 3.4.11]{NPS}). The last condition together with  the refinement equation (\ref{refeq}) and continuity of $\widehat{\varphi}$ at the origin imply that $\varphi$ forms an MRA (see, \cite[Theorem 1.2.14]{NPS}).
The refinement mask $m_0$ is also responsible for the stability of integer shifts of the refinable function $\varphi.$ We recall here necessary notions.  Let $g(\xi)$ be a $2\pi$-periodic function. If for $\alpha\in \mathbb{T}$ we have 
$g(\alpha)= g(\alpha+\pi)=0,$ then the pair $\{\alpha,\,\alpha+\pi\}$ is called a pair of symmetric roots of $g(\xi).$ A set of different complex numbers $\{b_1,\dots,b_n\}$ is called  cyclic if $b_{j+1}=b_j^2,$ $j=1,\dots,n$, and $b_n^2=b_1.$ In other words, $\displaystyle b_j={\rm e}^{i \beta_j}$, where  $\beta_j=2\pi m 2^{j-1}(2^n-1)^{-1}$ for some $m\in\mathbb{Z}.$ A cyclic set is called a cycle of a function $g$ if  $g(\beta_j+\pi)=0$ for all $j=1,\dots,n.$ The cycle $\{1\}$ is called trivial, all other cycles are called nontrivial.
 To provide the stability of the integer shifts we use the following statement \cite[Corollary 3.4.15]{NPS}: Integer shifts of a refinable function are stable if and only if the mask $m_0$ has neither nontrivial cycles nor a pair of  symmetric roots on $\mathbb{T}$.

To construct our wavelet frame we find a trigonometric polynomial $m_0$  approximating a continuous function. This trigonometric polynomial plays a role of a refinement mask of the UEP. It is impossible   to approximate an arbitrary continuous function $f$ by $m_0$ straightforward via  the distance $\|f-m_0\|_{\infty}$, since the system (\ref{matr}) implies the crucial  inequality   on the refinement function 
\begin{equation} 
\label{refmask}
\left|m_0(\xi)\right|^2+\left| m_0(\xi+\pi)\right|^2 \le 1.
\end{equation}
Therefore,  $|m_0|$ can not approximate a function that exceeds $1$.  That is why we introduce an auxiliary function $\tau$, which plays a role of a ``correction factor'' for the refinement mask. It turns out that for any continuous function $f$ such that $f(0)=1,$ and for any $\varepsilon>0$ there exists an integer $K_0 =K_0(f)$ and a mask $m_0$ such that 
$\|f - \tau^{-K_0}m_0\|<\varepsilon.$    

Suppose $m_0$ is a trigonometric polynomial, $m_0(0)=1$, and (\ref{refmask}) holds. 
 These conditions are sufficient for 
$m_0$ to be  a refinement mask of the UEP. Indeed, define $\displaystyle \widehat{\varphi}(\xi) = \prod_{j=1}^{\infty} m_0(\xi/2^j)$. Then $\varphi$ is the refinement function corresponding to the refinement mask $m_0(\xi)$: by the Mallat theorem \cite[Lemma 4.1.3]{NPS})  $\varphi \in L_2(\mathbb{R})$, and since $\widehat{\varphi}$ is an entire function of exponential  type, it follows that  $\widehat{\varphi}$ is continuous at zero and  
$\displaystyle \widehat{\varphi}(0) = m_0(0)=1$.

In \cite {CH} and \cite{Petukhov} it is proved that the inequality (\ref{refmask}) provides a frame with two generators, moreover, if $m_0$ is a trigonometric polynomial, it is possible to choose  $m_1$ and $m_2$ as the polynomials of the same degree.  We propose here one more method to construct wavelet masks. While we obtain three wavelet generators, not two, our formulas for wavelet masks are extremely simple and seems they could be interesting for applications.  

\begin{lem}
\label{wavemask}
Let $m_0$ be a $2\pi$-periodic trigonometric polynomial, $m_0(0)=1$, and  inequality (\ref{refmask}) holds. Then the functions 
 $$
m_1(\xi)={\rm e}^{ i \xi} \overline{m_0(\xi +\pi)},
\ \ 
m_2(\xi)=b(\xi)\cos (\xi/2), \ \ m_3(\xi)={\rm e}^{ i \xi} \overline{m_2(\xi +\pi)},
$$ 
where $|b(\xi)|^2=1-|m_0(\xi)|^2-|m_0(\xi+\pi)|^2,$  are trigonometric polynomials satisfying  system (\ref{matr}) with $q=3.$
\end{lem}

\textbf{Proof}. The function $b$ can be chosen as a trigonometric polynomial by the Riesz lemma. 
Therefore, the functions  $m_r$, $r=1,2,3,$ are trigonometric polynomials. We get 
$$
m_1(\xi)\overline{m_1(\xi+\pi)} = - m_0(\xi)\overline{m_0(\xi+\pi)},
$$  
$$
m_3(\xi)\overline{m_3(\xi+\pi)} = - m_2(\xi)\overline{m_2(\xi+\pi)}.
$$
Thus, the second equation is fulfilled. For the first equation we obtain
$$
\sum\limits_{r=0}^{3} |m_r(\xi)|^2 = 
 |m_0(\xi)|^2 + |m_0(\xi+\pi)|^2 +
|b(\xi)|^2 \cos^2 (\xi/2) + |b(\xi+\pi)|^2 \sin^2 (\xi/2).  
$$
 Since $|b|^2$ is a $\pi$-periodic function, it follows that the last expression is 
$$ 
 |m_0(\xi)|^2 + |m_0(\xi+\pi)|^2 +
|b(\xi)|^2\left(\cos^2 (\xi/2) + \sin^2 (\xi/2) \right) = |m_0(\xi)|^2 + |m_0(\xi+\pi)|^2 +
|b(\xi)|^2=1.  
$$
\hfill $\Box$ 

Now we are ready to prove the main result of the paper. 
 
\begin{teo}
Suppose $f\in C$, $f(0)=1,$  $\varepsilon >0$, $N_0\in \mathbb{N}$. Then there exists $K_0 = K_0(f)\in \mathbb{N}$ and a compactly supported Parseval wavelet frame with a refinement mask $m_0$ such that $\|f-\tau^{-K_0} m_0\|_{C}<\varepsilon,$ where   
$\displaystyle 
\tau(\xi):=1 -  \left(\int\limits_{0}^{\pi}\sin^{2N_0+1}t \, d t\right)^{-1} \int\limits_{0}^{\xi}\sin^{2N_0+1}t \, dt.
$ 
 The designed wavelet system has  vanishing moments of order $N_0$. The refinable function $\varphi$ has stable integer shifts.
\label{main}
\end{teo} 
 
\textbf{Proof}. 1.  
We approximate the function $f$ by a piecewise linear function with only finite number of roots, without nontrivial cycles and without pairs of symmetric roots. To this end, we find a piecewise linear function $f_1$ such that $\|f-f_1\|<\varepsilon/6.$ 
For example, $f_1$ interpolates $f$ at the equidistant points 
$\{\xi_k\}_{k=1}^{n},$ $-\pi=\xi_1<\xi_2< \dots<\xi_n=\pi.$  
Now we change the values of $f_1$ in the neighborhoods of the  segments, where $f_1\equiv 0$, if any. 
Let $f(\xi_k)=0$ for $k=i,\dots,j,$ and $f(\xi_{i-1}) \neq 0,$ $f(\xi_{j+1}) \neq 0.$ There are three possibilities. 

1) If $f_1(\xi)>0$ for $\xi\in (\xi_{i-1},\,\xi_i)$ and $f_1(\xi)>0$ for 
$\xi\in (\xi_{j},\,\xi_{j+1})$, then we define a function $f_2$ as
$f_2(\xi)=\min\{\varepsilon/12,\,f(\xi_{i-1}),\,f(\xi_{j+1})\}=:\gamma_1$ for $\xi\in[\xi^{\ast}_1,\,\xi^{\ast \ast}_1],$ where 
$$
\xi^{\ast}_1=\frac{\gamma_1 (\xi_{i-1}-\xi_i)}{f(\xi_{i-1})}+\xi_i \ \ \   \mbox{ and } \ \ \   \xi^{\ast \ast}_1=\frac{\gamma_1 (\xi_{j+1}-\xi_j)}{f(\xi_{j+1})}+\xi_j.
$$ 

2) If $f_1(\xi)<0$ for $\xi\in (\xi_{i-1},\,\xi_i)$ and $f_1(\xi)<0$ for 
$\xi\in (\xi_{j},\,\xi_{j+1})$, then analogously to 1) we define 
$f_2(\xi)=\max\{-\varepsilon/12,\,f(\xi_{i-1}),\,f(\xi_{j+1})\}=:\gamma_2$ for $\xi\in[\xi^{\ast}_2,\,\xi^{\ast \ast}_2],$ where 
$$
\xi^{\ast}_2=\frac{\gamma_2 (\xi_{i-1}-\xi_i)}{f(\xi_{i-1})}+\xi_i \ \ \   \mbox{ and } \ \ \   \xi^{\ast \ast}_2=\frac{\gamma_2 (\xi_{j+1}-\xi_j)}{f(\xi_{j+1})}+\xi_j.
$$ 

3) If the signs of $f_1$ are different on the intervals 
$(\xi_{i-1},\,\xi_i)$ and $(\xi_{j},\,\xi_{j+1})$, say, $f_1$ is positive on 
$(\xi_{i-1},\,\xi_i)$ and $f_1$ is negative on $(\xi_{j},\,\xi_{j+1})$, then we define $f_2$ as 
$$ \displaystyle
 f_2(\xi) = 
\left\{
\begin{array}{ll}
\frac{\gamma_3-f(\xi_{i-1})}{\xi_{i}-\xi_{i-1}} (\xi - \xi_{i}) + \gamma_3,  & \xi\in[\xi_{i-1},\,\xi_{i}], \\
\frac{\gamma_3-\gamma_4}{\xi_{j}-\xi_{i}} (\xi - \xi_{i}) + \gamma_3,  & \xi\in[\xi_{i},\,\xi_{j}], \\
\frac{f(\xi_{j+1})-\gamma_4}{\xi_{j+1}-\xi_{j}} (\xi - \xi_{j}) + \gamma_4,  & \xi\in[\xi_{j},\,\xi_{j+1}], \\
\end{array}
\right.
$$
where $\gamma_3:=\min\{\varepsilon/12,\,f(\xi_{i-1})\},$ $\gamma_4:=\max\{-\varepsilon/12,\,f(\xi_{j+1})\}.$
On the remaining part of $\mathbb{T}$ the functions $f_1$ and $f_2$ coincides.

Roughly speaking, 
if $f_1$ has the same sign to the left of the point $\xi_{i-1}$ and to the right of the point $\xi_{j+1}$, then we ``shift up'' (if the sign is ``plus'') or ``shift down'' (if the sign is ``minus'') the plot of $f_1$ to $\gamma_1$ or $\gamma_2$. If the signs of $f_1$ are different, 
then we ``rotate'' the plot of $f_1$. As a result, we get the piecewise linear function $f_2$ that has only finite number of roots and $\|f-f_2\|_C<\varepsilon /3.$ We keep notation $\xi_k$ for the first coordinates of  nodes for $f_2.$ 

Finally, we remove pairs of symmetric roots and cycles of $f_2$ if any. Let $\xi_0$ be one of the symmetric roots or one of the roots from a cycle and $\xi_0\in[\xi_k,\,\xi_{k+1}).$ 
If $\xi_0=\xi_k$, then we shift the node of interpolation $(\xi_k,\,0)$ slightly keeping the order of the first coordinates of nodes. In other words,  for the new node $(\tilde{\xi}_k,\,0)$ we get $\tilde{\xi}_k \in (\xi_{k-1},\,\xi_{k+1}).$  If $\xi_0 \neq \xi_k$, then we shift  the root replacing one link of the polyline with two. More precisely, we replace the link $AB$, where $A=(\xi_k,\,f(\xi_k)),$ $B=(\xi_{k+1},\,f(\xi_{k+1})),$ with two links $AC$ and $CB$, where $C=(\xi_0,\,h)$, $|h|<\varepsilon/6.$ Thus, in both cases  we get a new root $\xi'_0$.  As a result, we obtain a new piecewise linear function
 $f_3$ such that $\|f-f_3\|_C< \varepsilon/2$. This function has only finite number of roots. It has  neither nontrivial cycles nor pairs of symmetric roots on $\mathbb{T}$.         

2. We find a  trigonometric polynomial $T$ such that 
$\|f_3-T\|_C<\varepsilon_1$,  $T(0)=1,$ $T^{(j)}(0)=0,$ $j=1,\dots,J,$ and $T$ has neither nontrivial cycles nor pairs of symmetric roots on $\mathbb{T}$. 
To this end,  let $T_1$ be a  trigonometric polynomial  approximating the function $f_3$, so $\|f_3-T_1\|_C<\varepsilon_1/2$.    Let $T_1^{(j)}(0)=\gamma_j,$ $j=1,\dots,J.$ 
We need to find a polynomial 
$$ 
T_2(\xi)= \sum_{j=1}^{\lfloor J/2 \rfloor}\alpha_j \cos \,  N_j \xi +
\sum_{j=1}^{\lfloor (J+1)/2 \rfloor}\beta_j \sin \,  K_j \xi
$$
such that $\|T_2\|_C<\varepsilon_1/2$ and $T_2^{(j)}(0)=-\gamma_j,$ $j=1,\dots,J.$ 
Then $\|f_3-(T_1+T_2)\|_C<\varepsilon_1$ and $(T_1+T_2)^{(j)}(0)=0,$ $j=1,\dots,J.$
To find  $T_2$ we 
interpret  the conditions $T_2^{(j)}(0)=-\gamma_j$ as two linear systems  in the variables $\alpha_j,$ $\beta_j$:
$$
\sum_{i=1}^{\lfloor J/2 \rfloor} \alpha_i N_i^{2k}=(-1)^k \gamma_{2k},
\ \ k=1,\dots,\lfloor J/2 \rfloor,
$$
$$
\sum_{i=1}^{\lfloor (J+1)/2 \rfloor} \beta_i K_i^{2k-1}=(-1)^{k-1} \gamma_{2k-1}, \ \ k=1,\dots,\lfloor (J+1)/2 \rfloor.
$$
Let us consider the first system. For  the second system, we can proceed the same way. The determinant of the first system is the determinant of a square Vandermonde matrix and it is not equal to zero iff $N_i \neq N_j$ as  $i \neq j,$ thus, there exists a unique solution. It has the form 
$$
\sum_{k=1}^{\lfloor J/2 \rfloor} (-1)^k \gamma_{2k} n_{i,k},
$$
where  the element of the inverse matrix of the system, $n_{i,k}$,  
is equal to a rational function in variables $N_j,$ $j=1,\dots,\lfloor J/2 \rfloor$, the denominator of the function is a homogeneous polynomial of the order $\lfloor J/2 \rfloor \lfloor J/2+1 \rfloor$. And the main point is that the order of the denominator is strictly greater then the order of the nominator. Therefore, we can choose $N_j$ large enough to provide the inequalities 
$|\alpha_i|<\varepsilon_1/(2J).$ Analogously, we can choose $K_j$ such that 
$|\beta_i|<\varepsilon_1/(2J).$ 
  Then $\|T_2\|_C<\sum_{i=1}^{\lfloor J/2 \rfloor} |\alpha_i|+ 
\sum_{i=1}^{\lfloor (J+1)/2 \rfloor} |\beta_i|<\varepsilon_1/2.$
 We set 
$T(\xi):=(T_1(\xi)+T_{2}(\xi))/(T_1(0)+T_{2}(0)).$ 

Now we claim that for a small enough $\varepsilon_1$ the polynomial $T$ has neither nontrivial cycles nor pairs of symmetric roots on $\mathbb{T}$. Indeed, suppose $f_3$ is defined explicitly as $f_3(\xi)=a_k \xi+b_k$ for $\xi \in [\xi_k,\xi_{k+1}],$ $k=1,\dots,n-1,$ $\xi_1=-\pi,$ $\xi_{n}=\pi.$ 
 We recall that  $\xi_0$ is a root and  we replace it to the new root $\xi'_0$ to remove pairs of symmetric roots and cycles of the function $f_2$. We denote by  $\alpha$ the minimum of $|\xi_0-\xi'_0|/2$ over all the pairs of the old and new roots $(\xi_0,\,\xi'_0)$  and  $a:=\min\{|a_k|:a_k \neq 0\}.$  Suppose the function $f_3$ has a root $\xi^0$ on the segment $[\xi_k,\xi_{k+1}]$. If $\|f_3-T\|_C<\varepsilon_1$, then 
$a_k \xi+b_k-\varepsilon_1 \leq T(\xi) \leq a_k \xi+b_k+\varepsilon_1$. In other words,  the plot of the polynomial $T$ lies inside the parallelogram bounded by the lines $y(\xi)=a_k \xi+b_k \pm \varepsilon_1$, $\xi=\xi_k$, $\xi=\xi_{k+1}.$ 
Therefore, the roots of $T$ can only be in the neighborhood of $\xi^0,$ namely in the interval of the length $2\varepsilon_1/a_k$. If we provide the inequality 
$\varepsilon_1/a_k<\alpha$ for all the segments  $[\xi_k,\xi_{k+1}]$ that contains roots of $f_3$, then  it do means that $T$ has neither nontrivial cycles nor pairs of symmetric roots on $\mathbb{T}$. To provide the inequality it is sufficient   to choose $\varepsilon_1< a \alpha.$

3. To design a refinement mask from the polynomial $T$ we  consider the $2\pi$-periodic function
$$ 
\tau(\xi):=1 - \frac{\int\limits_{0}^{\xi}\sin^{2N_0+1}t \, dt}{\int\limits_{0}^{\pi}\sin^{2N_0+1}t \, dt},
$$ 
 and 
  find $K_0 \in \mathbb{N}$ such that for all $\xi$
\begin{equation}\label{less1}
\left(\tau (\xi)\right)^{2K_0}\left|T(\xi)\right|^2 + 
\left(\tau(\xi+\pi)\right)^{2K_0}\left|T(\xi+\pi)\right|^2 \le 1.
\end{equation} 
Since the function $A(\xi):=\left(\tau (\xi)\right)^{2K_0}\left|T(\xi)\right|^2 + 
\left(\tau(\xi+\pi)\right)^{2K_0}\left|T(\xi+\pi)\right|^2$ is $\pi$-periodic,  it is sufficient to check (\ref{less1}) for $\xi \in [-\pi/2,\,\pi/2).$ For the function $\tau$ we get 
$$\displaystyle \tau'(\xi)=-\left(\int\limits_{0}^{\pi}\sin^{2N_0+1}t \, dt\right)^{-1} \sin^{2N_0+1}\xi = -\left(\int\limits_{0}^{\pi}\sin^{2N_0+1}t \, dt\right)^{-1} \xi^{2N_0+1} + o(\xi^{2N_0+2})$$ and $\tau(0)=1$. So,
\begin{equation}
\label{est1}
\tau(\xi) = 1 -c \xi^{2N_0+2} + o(\xi^{2N_0+2}),
\end{equation}
  where $\displaystyle c:=\left((2N_0+2)\int\limits_{0}^{\pi}\sin^{2N_0+1}t \, dt\right)^{-1}>0.$ 
	Analogously, 
	\begin{equation}
	\label{zero_pi}
	\tau(\xi+\pi) = c\, \xi^{2N_0+2} + o(\xi^{2N_0+2}).
	\end{equation}
	It follows from item 1. that  
\begin{equation}
\label{est2}
T(\xi)=1+a\,\xi^{J+1}+o(\xi^{J+1}).
\end{equation}
In the sequel we choose $J$ to satisfy the inequality  $J>2N_0+1.$ Thus, for the function $A_K$, we have 
$$
A_K(\xi) = \left(1 -c \,\xi^{2N_0+2} + o(\xi^{2N_0+2})\right)^{2K} \left(1+a\,\xi^{J+1}+o(\xi^{J+1})\right)^2
$$
$$ + \left( c\, \xi^{2N_0+2} + o(\xi^{2N_0+2})\right)^{2K} \left|T(\xi+\pi)\right|^2 
=1-2Kc\, \xi^{2N_0+2} +o(\xi^{2N_0+2}).
$$
Therefore, $A_K(0)=1$ is a local maximum of $A_K.$ Thus, there exists a neighborhood $U$ of the point $0$ such that $A_1(\xi)\le 1$ for $\xi \in U.$ It remains to note that $A_{K_1}(\xi)\leq A_{K_2}(\xi)$ for $K_1\ge K_2$ and $A_K(\xi)$   converges uniformly to $0$ on $[-\pi/2,\,\pi/2)\setminus U$ as $K\to \infty.$ Therefore, there exists an integer $K_0$ such that $A_K(\xi)\le 1$ on $[-\pi/2,\,\pi/2)\setminus U$ for $K\ge K_0$. 

We define the $2\pi$-periodic  trigonometric polynomial $m_0(\xi):= \left(\tau(\xi)\right)^{K_0} T(\xi).$  Since $m_0(0)=1$ and (\ref{less1}) holds, it follows that 
the function $\displaystyle \widehat{\varphi}(\xi):= \prod_{j=1}^{\infty}m_0(\xi 2^{-j})\in L_2(\mathbb{R})$ is the Fourier transform of the refinable function corresponding to the refinement mask $m_0.$ Defining the wavelet masks $m_r,$ $r=1,\dots,q$ via the UEP as it is done in   \cite{CH} or \cite{Petukhov}, or Lemma \ref{wavemask}, we obtain the wavelet generators $\psi_r,$ $r=1,\dots,q$.

4. The last we need to provide is the vanishing moments of the wavelet generators. 
 Since
 $
\widehat{\psi}_r(2\xi)= m_r(\xi)\widehat{\varphi}(\xi),
$ and $\widehat{\varphi}(0)=1$,
we need to check 
$m_r^{(j)}(0)=0 $  for   $j=0,\dots, N_0.$ And taking into account the equality 
$
\sum\limits_{r=1}^{q} |m_r(\xi)|^2 = 1-|m_0(\xi)|^2,
$
it is necessary and  sufficient to provide 
$\left(1-|m_0(\xi)|^2\right)^{(j)}(0)=0 $  for   $j=0,\dots, 2N_0.$
It remains to use the estimates (\ref{est1}), (\ref{est2}) for the function $A_K$ from the item 2. So we get 
$$
1-|m_0(\xi)|^2 = 1- \left(\tau(\xi)\right)^{2K_0} (T(\xi))^2
$$
$$
=1-\left(1 -c\, \xi^{2N_0+2} + o(\xi^{2N_0+2})\right)^{2K_0} \left(1+a\, \xi^{J+1}+o(\xi^{J+1})\right)^2
= 2 K_0 c\, \xi^{2N_0+2} +o(\xi^{2N_0+2}).
$$
 \hfill $\Box$

\begin{rem}
We recall that for the compactly supported refinable function with stable integer shifts the approximation order $N_{\varphi}$ is defined by the order of the zero of the mask  $m_0$   at the point $\xi=\pi$. Namely, if the order of the zero is $\ge n$, then the approximation order is $\ge n$ (see \cite[Theorems  3.3.2 and 3.4.16]{NPS}). In our case the order of the zero of $m_0$ at $\xi=\pi$ is $\ge$  the order of the zero of $\tau^{K_0}$ at $\xi=\pi$. According to (\ref{zero_pi}) the last is equal to $K_0(2N_0+2).$        
\end{rem}

 The work  is supported by the
Russian Science Foundation (grant 18-11-00055).

\end{document}